\documentclass[12pt,leqno]{article}
\usepackage{latexsym}
\usepackage{amsmath}
\usepackage{amssymb}

\newcommand{\bi}{\begin{itemize}}
\newcommand{\ei}{\end{itemize}}
\newcommand{\beq}{\begin{equation}}
\newcommand{\eeq}{\end{equation}}
\newcommand{\abso}[1]{\left| #1 \right|}

\newcommand{\fin}{\begin{flushright} $\Box $ \end{flushright}}
\newcommand{\vsp}{\vspace{10mm}}

\newcommand{\matR}{\mathbb{R}}

\newtheorem{thm}{THEOREM}
\newtheorem{dfn}{DEFINITION}

\newtheorem{lm}{LEMMA}
\newtheorem{cor}{COROLLARY}

\newcounter {nrq} 

\begin{document}

\title{\bf\sc L\'{e}vy Approximation of Impulsive Recurrent Process
with Semi-Markov Switching}

\author{{\sc V. S. Koroliuk}$^1$, {\sc N. Limnios}$^2$ and {\sc I.V. Samoilenko}$^1$\\
$^1$Institute of Mathematics,\\ Ukrainian National Academy of Science, Kiev, Ukraine\\
$^2$Laboratoire de Math\'ematiques Appliqu\'ees,\\ Universit\'e de Technologie de Compi\`egne, France}


\maketitle

\baselineskip 6 mm

\vsp

\hrule
\begin{abstract}
In this paper, the weak convergence of impulsive recurrent process
with semi-Markov switching in the scheme of L\'{e}vy approximation
is proved. Singular perturbation problem for the compensating
operator of the extended Markov renewal process is used to prove the
relative compactness.
\end{abstract}


\vsp
\hrule
\section{Introduction}
L\'{e}vy approximation is still an active area of research in
several theoretical and applied directions. Since L\'{e}vy processes
are now standard, L\'{e}vy approximation is quite useful for
analyzing complex systems (see, e.g. \cite{ber, sat}). Moreover they
are involved in many applications, e.g., risk theory, finance,
queueing, physics, etc. For a background on L\'{e}vy process see,
e.g. \cite{ber, sat, gisk}.

In particular in \cite{korlim} it has been studied the following
impulsive process as partial sums in a series scheme
\begin{eqnarray}\label{1adf1}
\xi^{\varepsilon}(t)=\xi_0^{\varepsilon}+\sum_{k=1}^{\nu(t)}\alpha^{\varepsilon}_k(x^{\varepsilon}_{k-1}),\quad
t\ge 0,
\end{eqnarray}
the random variables $\alpha_k^{\varepsilon} (x), k \geq 1$ are
supposed to be independent and perturbed by the jump Markov process
$x(t), t\ge 0$.

We propose to study generalization of the problem (\ref{1adf1}):
\begin{eqnarray}\label{1bdf1}
\xi^{\varepsilon}(t)=\xi_0^{\varepsilon}+\sum_{k=1}^{\nu(t)}\alpha^{\varepsilon}_k(\xi^{\varepsilon}_{k-1},x^{\varepsilon}_{k-1}),\quad
t\ge 0.
\end{eqnarray}
Here the random variables $\alpha_k^{\varepsilon} (u, x), k \geq 1$
depend on the process $\xi^{\varepsilon}(t)$.

We propose to study convergence of (\ref{1bdf1}) using a combination
of two methods. The one, based on semimartingales theory, is
combined with a singular perturbation problem for the compensating
operator of the extended Markov renewal process. So, the method
includes two steps.

In the first step we prove the relative compactness of the
semimartingales representation of the family $\xi^\varepsilon$,
$\varepsilon>0$, by proving the following two facts \cite{ethier}:
$$\lim\limits_{c\to \infty}\sup\limits_{\varepsilon\leq\varepsilon_0}
\mathbf{P}\{\sup\limits_{t\leq T}|\xi^{\varepsilon}(t)|>c\}=0,$$
known as the compact containment condition,
and $$\mathbf{E}|\xi^{\varepsilon}(t)-\xi^{\varepsilon}(s)|^2\le k |t-s|,$$
for some positive constant $k$.

In the second step we prove convergence of the extended Markov
renewal process $\xi_n^{\varepsilon}, x_n^{\varepsilon},
\tau_n^{\varepsilon}, n\geq0$ by using singular perturbation
technique as presented in \cite{korlim}.

Finally, we apply Theorem 6.3 from \cite{korlim}.

The paper is organized as follows. In Section 2 we present the
time-scaled impulsive process (\ref{1bdf1}) and the switching
semi-Markov process. In the same section we present the main results
of L\'{e}vy approximation. In Section 3 we present the proof of the
theorem.

\vsp \hrule
\section{Main results}
Let us consider the space $\matR^d$ endowed with a norm $\abso{\cdot}$ ($d\ge 1$),
and $(E,\mathcal{E})$, a {\it standard phase space}, (i.e.,
$E$ is a Polish space and $\mathcal{E}$ its Borel $\sigma$-algebra).
For a vector $v\in \matR^d$ and a matrix $c\in \matR^{d\times d}$ ,
$v^*$ and $c^*$ denote their transpose respectively.
Let $C_3(\matR^d)$ be a measure-determining class of real-valued bounded functions, such that
$g(u)/\abso{u}^2 \to 0$, as $\abso{u}\to 0$ for $g\in C_3(\matR^d)$ (see \cite{jacod1,korlim}).

The impulsive processes $\xi^{\varepsilon}(t), t\geq 0,
\varepsilon>0$ on $\mathbb{R}^d$ in the series scheme with small
series parameter $\varepsilon\to 0$, $(\varepsilon>0)$ are defined
by the sum (\cite[Section 9.2.1]{korlim})
\begin{eqnarray}\label{1adf2}
    \xi^{\varepsilon}(t)=\xi_0^{\varepsilon}+\sum_{k=1}^{\nu(t/\varepsilon^2)}\alpha^{\varepsilon}_{k}(\xi^{\varepsilon}_{k-1},x^{\varepsilon}_{k-1}),\quad
t\ge 0.
\end{eqnarray}

For any $\varepsilon > 0$, and any sequence $z_k, k \geq 0$, of
elements of $\mathbb{R}^d\times E$, the random variables
$\alpha_k^{\varepsilon} (z_{k-1}), k \geq 1$ are supposed to be
independent. Let us denote by $G_{u,x}^{\varepsilon}$ the
distribution function of $\alpha_k^{\varepsilon} (x)$, that is,
$$G_{u,x}^{\varepsilon}(dv) := P(\alpha_k^{\varepsilon} (u,x) \in
dv), k \geq 0, \varepsilon
> 0, x \in E, u\in \mathbb{R}^d.$$ It is worth noticing that the coupled process $\xi^\varepsilon(t),
x^{\varepsilon}(t), t \geq 0$, is a Markov additive process (see,
e.g., \cite[Section 2.5]{korlim}).

We make natural assumptions for the counting process $\nu(t)$,
namely:
\begin{eqnarray}\label{1con1}
\int_0^t\mathbf{E}[\varphi(s)d\nu(s)]<l_1\int_0^t\mathbf{E}(\varphi(s))ds
\end{eqnarray}
for any nonnegative, increasing $\varphi(s)$ and $l_1>0.$

The switching semi-Markov process ${x}(t), t\ge 0$ on the standard
phase space $(E,\mathcal{E})$, is defined by the semi-Markov kernel
$$Q(x,B,t) = P(x,B)F_x(t), x \in E,B \in\mathcal{E}, t \geq 0, $$ which defines
the associated Markov renewal process $ x_n, \tau_n, n \geq 0 $:
$$Q(x,B, t) = P(x_{n+1} \in B, \theta_{n+1} \leq t | x_n = x) = P(x_{n+1} \in B | x_n =
x)P(\theta_{n+1}\leq t | x_n = x).$$

Finally we should denote $\xi^{\varepsilon}_n$ in (3):
$$\xi^{\varepsilon}_n:=\xi(\varepsilon^2\tau_n)=\xi_0^{\varepsilon}+\sum_{k=1}^n\alpha_k^{\varepsilon}(\xi^{\varepsilon}_{k-1},x^{\varepsilon}_{k-1}).$$

The L\'{e}vy approximation of Markov impulsive process (\ref{1adf2})
is considered under the following conditions.

\begin{description}
\item[C1:] The semi-Markov process ${x}(t), t \geq 0$ is uniformly
ergodic with the stationary distribution
$$\pi(dx)q(x)=q\rho(dx), q(x):=1/m(x), q:=1/m,$$
$$m(x):=\mathbb{E}\theta_x=\int_0^{\infty}\overline{F}_x(t)dt, m:=\int_E\rho(dx)m(x),$$
$$\rho(B)=\int_E\rho(dx)P(x,B), \rho(E)=1.$$

\item[C2:] {\it L\'{e}vy approximation}. The family of impulsive processes $\xi^{\varepsilon}(t), t\geq 0$ satisfies the
L\'{e}vy approximation conditions \cite[Section 9.2]{korlim}.
\begin{description}
\item[L1:] Initial value condition
$$\sup\limits_{\varepsilon>0} E|\xi_0^{\varepsilon}|\leq C <
\infty$$ and $$\xi_0^{\varepsilon}\Rightarrow\xi_0.$$

\item[L2:]Approximation of the mean values:
$$a^{\varepsilon}(u;x) = \int_{\mathbb{R}^d} vG^{\varepsilon}_{u,x}(dv)
= \varepsilon a_1(u;x)+\varepsilon^2 [a(u;x) +\theta_a^{\varepsilon}
(u;x)],$$ and
$$c^{\varepsilon}(u;x) = \int_{\mathbb{R}^d}
vv^*G^{\varepsilon}_{u,x}(dv) = \varepsilon^2 [c(u;x) +
\theta_c^{\varepsilon} (u;x)],$$ where functions $a_1, a$ and $c$
are bounded.

\item[L3:] Poisson approximation condition for intensity kernel (see \cite{jacod1})
$$G_g^{\varepsilon}(u;x) = \int_{\mathbb{R}^d} g(v)G^{\varepsilon}_{u,x}(dv)
= \varepsilon^2[G_g(u;x) + \theta^{\varepsilon}_g(u;x)]$$ for all $g
\in C_3(\mathbb{R}^d)$, and the kernel $G_g(u;x)$ is bounded for all
$g \in C_3(\mathbb{R}^d)$, that is,
$$|G_g(u;x)| \leq G_g \quad \hbox{(a constant depending on $g$)}.$$

Here \begin{eqnarray}\label{1gg2}
   G_{g}(u;x) =\int_{\mathbb{R}^d} g(v)G_{u,x}(dv),\quad g \in C_3(\mathbb{R}^d).
\end{eqnarray}

The above negligible terms
$\theta_a^\varepsilon,\theta_c^\varepsilon, \theta_g^\varepsilon$
satisfy the condition $$\sup\limits_{x\in E}
|\theta_{\cdot}^{\varepsilon}(u;x)|\to 0,\quad  \varepsilon\to 0.$$

\item[L4:] {\it Balance condition}.
$$\int_E\rho(dx)a_1(u;x)=0.$$
\end{description}

In addition the following conditions are used:
\item[C3:] {\it Uniform square-integrability}:
$$\lim\limits_{c\to\infty}\sup\limits_{x\in E} \int_{|v|>c} vv^*G_{u,x}(dv) = 0.$$

\item[C4:] {\it Linear growth}: there exists a positive constant $L$ such that
$$|a(u;x)|\leq L(1+|u|),\quad\hbox{and}\quad |c(u;x)|\leq L(1+\abso{u}^2),$$
and for any real-valued non-negative function $f(v), v\in
\mathbb{R}^d$, such that $\int_{\mathbb{R}^d\setminus
\{0\}}(1+f(v))\abso{v}^2dv<\infty,$ we have
$$|G_{u,x} (v)|\leq Lf(v)(1+\abso{u}).$$
\end{description}

\vsp
The main result of our work is the following.

\begin{thm} Under conditions $\mathbf{C1-C4}$
the weak convergence
$$\xi^{\varepsilon}(t)\Rightarrow \xi^0(t),\quad \varepsilon \to 0$$ takes
place.

The limit process $\xi^0(t), t\geq0$ is a L\'{e}vy process defined
by the generator $\mathbf{L}$ as follows
\begin{eqnarray}\label{1limgen}
\mathbf{L}\varphi(u)=(\widehat{a}(u)-\widehat{a}_0(u))\varphi'(u)+
\frac{1}{2}\sigma^2(u)\varphi''(u) + \lambda(u)\int_{\mathbb{R}^d}
[\varphi(u + v)-\varphi(u)]G_{u}^0(dv),
\end{eqnarray}
where:
$$\widehat{a}(u)=q\int_E\rho(dx)a(u;x), \widehat{a}_0(u)=\int_EvG_u(dv), G_u(dv)=q\int_E\rho(dx)G_{u,x}(dv),$$
$$\widehat{a^2_1}(u)=q\int_E\rho(dx)a^2_1(u;x),\hskip 5mm \widetilde{a}_1(u;x):=q(x)\int_EP(x,dy)a_1(u;x), c_0(u;x)=\int_Evv^*G_{u,x}(dv)$$
$$\sigma^2(u)=2\int_E\pi(dx)\{\widetilde{a}_1(u;x)\widetilde{R}_0\widetilde{a}_1^*(u;x)+
\frac{1}{2}[c(u;x)-c_0(u;x)]\}-\widehat{a^2_1}(u), \hskip 5mm
\sigma^2(u)\geq 0$$
\hskip40mm $\lambda(u)=G_u(\mathbb{R}^d),$ \hskip5mm
$G_{u}^0(dv)=G_u(dv)/\lambda(u),$

here $\widetilde{R}_0$ is the potential operator of embedded Markov
chain.
\end{thm}

\textbf{Remark 1.} The limit L\'{e}vy process consists of three
parts: deterministic drift, diffusion part and Poisson part.

There are some possible cases:

\begin{description}
\item[1).] If $\widehat{b}(u)-\widehat{b}_0(u)=0$ then the limit
process does not have deterministic drift.

\item[2).] If $\sigma^2(u)=0$ then the limit process does not have
diffusion part. As a variant of this case we note that if
$c(u;x)=c_0(u;x)$ then also $b_1(u;x)=0$ and we obtain the
conditions of Poisson approximation after re-normation
$\varepsilon^2=\widetilde{\varepsilon}$ (see, for example Chapter 7
in \cite{korlim}).
\end{description}

\textbf{Remark 2.} In the work \cite{korlim} (Theorem 9.3) an
analogical result was obtained for impulsive process with Markov
switching. If we study an ordinary impulsive process without
switching, we should obtain
$\sigma^2=E(\alpha_k^{\varepsilon})^2-(E(\alpha_k^{\varepsilon}))^2=(c-c_0)-a_1^2$.
This result correlates with the similar results from \cite{jacod1}.
In case of our Theorem this may be easily shown, but in
\cite{korlim} (Theorem 9.3) it is not obvious.

The difference is that we used $\widetilde{R}_0$ -- the potential
operator of embedded Markov chain instead of $R_0$ -- the potential
operator of Markov process. Due to this, our result obviously
correlates with other well-known result.

\textbf{Remark 3.} Asymptotic of the second moment in the condition
\textbf{L1} contains second modified characteristics $c(u;x)$ (see
correlation 4.2 at page 555 in \cite{jacod1}). This characteristics
in limit contains both second moment of Poisson part and dispersion
of diffusion part, namely $c=c_0+\sigma^2.$

\vsp
\hrule
\section{Proof of Theorem 1}
The proof of Theorem 1 is based on the semimartingale representation
of the impulsive process (\ref{1adf2}).

We split the proof of Theorem 1 in the following two steps.

\noindent {\sc Step 1}. In this step we
establish the relative compactness of the family of processes
$\xi^{\varepsilon}(t), t\geq 0, \varepsilon>0$ by using the approach
developed in \cite{lip3}. Let us remind that the space of all probability
measures defined on the standard space $(E,{\cal E})$ is also a Polish space;
so the relative compactness and tightness are equivalent.

First we need the following lemma.

\begin{lm} Under assumption $\mathbf{C4}$ there exists a
constant $k>0$, independent of $\varepsilon$ and dependent on $T$, such that
$$\mathbf{E}\sup\limits_{t\leq T}|\xi^{\varepsilon}(t)|^2\leq k_T.$$
\end{lm}

\begin{cor}  Under assumption $\mathbf{C4}$, the following
compact containment condition (CCC) holds:
$$\lim\limits_{c\to \infty}\sup\limits_{\varepsilon\leq\varepsilon_0}
\mathbf{P}\{\sup\limits_{t\leq T}|\xi^{\varepsilon}(t)|>c\}=0.$$
\end{cor}
\noindent{\it Proof}: The proof of this corollary follows from Kolmogorov's inequality.\fin

\vsp \noindent{\it Proof of Lemma 1}: (following \cite{lip3}). The
impulsive process (\ref{1adf2}) has the following semimartingale
representation
\begin{eqnarray}\label{1smdecomp}
\xi^{\varepsilon}(t)=u+B_t^{\varepsilon}+M_t^{\varepsilon},
\end{eqnarray}
where $u= \xi^\varepsilon_0$; $B_t^{\varepsilon}$ is the predictable
drift
$$B_t^{\varepsilon}=\sum_{k=1}^{\nu(t/\varepsilon^2)}a^\varepsilon(\xi^{\varepsilon}_{k-1},{x}^{\varepsilon}_{k-1})
=A_1^\varepsilon(t)+A^\varepsilon(t)+\theta^\varepsilon_a(t),$$
where
$$A^{\varepsilon}_1(t):=\varepsilon\sum_{k=1}^{\nu(t/\varepsilon^2)}a_{1}
(\xi_{k-1}^{\varepsilon}, x_{k-1}^{\varepsilon}),
A^{\varepsilon}(t):= \varepsilon^2\sum_{k=1}^{\nu(t/\varepsilon^2)}a
(\xi_{k-1}^{\varepsilon} ,x_{k-1}^{\varepsilon}).$$

\begin{eqnarray}\label{1qch}
\langle M^{\varepsilon}\rangle_t=\varepsilon^2\sum_{k=1}^{\nu
(t/\varepsilon^2)}\int_{\mathbb{R}^d\setminus\{0\}}vv^*G(\xi^{\varepsilon}_{k-1},dv;{x}_{k-1}^{\varepsilon})+
\theta^{\varepsilon}_c(t)=\\ \nonumber \varepsilon^2\sum_{k=1}^{\nu
(t/\varepsilon^2)}
c(\xi^{\varepsilon}_{k-1};{x}_{k-1}^{\varepsilon})
+\theta^{\varepsilon}_c(t),
\end{eqnarray} and for every finite $T>0$
$$\sup\limits_{0\leq t\leq T} |\theta^\varepsilon_{\cdot}(t)|\rightarrow 0, \varepsilon\rightarrow 0.$$

To verify compactness of the process $\xi^{\varepsilon}(t)$ we split
it at two parts.

The first part of order $\varepsilon$
$$A_1^\varepsilon(t)=\varepsilon\sum_{k=1}^{\nu(t/\varepsilon^2)}
a_1(\xi^{\varepsilon}_{k-1};{x}^{\varepsilon}_{k-1}),$$ can be
characterized by the compensating operator
$$\mathbf{L}^{\varepsilon}\varphi(u;x)=\varepsilon^{-2}q(x)[\mathbf{A}_1^{\varepsilon}(x)P-I]\varphi(u;x),$$ where
$\mathbf{A}_1^{\varepsilon}(x)\varphi(u)=\varphi(u+\varepsilon
a_1(u;x))=\varepsilon a_1(u;x)\varphi'(u)+\varepsilon
\theta^{\varepsilon}\varphi(u).$ After simple calculations we may
rewrite the operator:
$$\mathbf{L}^{\varepsilon}=\varepsilon^{-2}\mathbf{Q}+\varepsilon^{-1}\mathbf{A}_1(x)P+\theta^{\varepsilon},$$
here $\mathbf{A}_1(x)\varphi(u)=\varepsilon a_1(u;x)\varphi'(u).$

Corresponding martingale characterization is the following
$$\mu_{n+1}^{\varepsilon}=\varphi(A_{1,n+1}^{\varepsilon},x_{n+1}^{\varepsilon})-\varphi(A_{1,0}^{\varepsilon},x_0^{\varepsilon})-
\sum_{m=0}^{\nu_n}\theta^{\varepsilon}_{m+1}\mathbf{L}^{\varepsilon}\varphi(A_{1,m}^{\varepsilon},x_m^{\varepsilon}).$$

Using the results from \cite{korlim}, Section 1 we obtain the last
martingale in the form
$$\widetilde{\mu}_t^{\varepsilon}=\varphi^{\varepsilon}(A_1^\varepsilon(t),x^{\varepsilon}_t)+
\varphi^{\varepsilon}(A_1^\varepsilon(0),x^{\varepsilon}_0)-
\int_0^t\mathbf{L}^{\varepsilon}\varphi^{\varepsilon}(A_1^\varepsilon(s),x^{\varepsilon}_s)ds,$$
where $x^{\varepsilon}_t:=x(t/\varepsilon^2).$

Thus (see, for example Theorem 1.2 in \cite{korlim}), it has
quadratic characteristic
$$<\widetilde{\mu}^{\varepsilon}>_t=\int_0^t\left[\mathbf{L}^{\varepsilon}(\varphi^{\varepsilon}
(A_1^\varepsilon(s),x^{\varepsilon}_t))^2-2\varphi^{\varepsilon}(A_1^\varepsilon(s),x^{\varepsilon}_s)
\mathbf{L}^{\varepsilon}\varphi^{\varepsilon}(A_1^\varepsilon(s),x^{\varepsilon}_s)\right]ds.$$

Applying the operator
$\mathbf{L}^{\varepsilon}=\varepsilon^{-2}\mathbf{Q}+\varepsilon^{-1}\mathbf{A}_1(x)P+\theta^{\varepsilon}$
to test-function
$\varphi^{\varepsilon}=\varphi+\varepsilon\varphi_1$ we obtain the
integrand of the view
$$Q\varphi_1^2-2\varphi_1Q\varphi_1+\theta^{\varepsilon}\varphi^{\varepsilon}.$$

Thus the integrand is limited. The boundedness of the quadratic
characteristic provides $\widetilde{\mu}_t^{\varepsilon}$ is
compact. Thus, $\varphi(A_1^\varepsilon(t))$ is compact too and
bounded uniformly by $\varepsilon$. By the results from
\cite{ethier} we obtain compactness of $A_1^\varepsilon(t)$, because
the test-function $\varphi(u)$ belongs to the measure-determining
class.

Now we should study the second part of order $\varepsilon^2$.

For a process $y(t), t\ge 0$, let us define the process
$y^\dag(t)=\sup\limits_{s\leq t}|y(s)|,$ then from (\ref{1smdecomp})
we have
\begin{eqnarray}\label{1eq4}
((\xi^{\varepsilon}(t))^\dag)^2\le
4[u^2+((A^{\varepsilon}(t))^\dag)^2+((M^{\varepsilon}_t)^\dag)^2].
\end{eqnarray}

Now we may apply the result of Section 2.3 \cite{korlim}, namely
$$\sum_{k=1}^{\nu(t)} a(\xi^\varepsilon_{k-1}, x^\varepsilon_{k-1})=\int_0^ta(\xi^{\varepsilon}(s),x^{\varepsilon}(s))d\nu(s).$$

Condition $\mathbf{C4}$  implies that for sufficiently large
$\varepsilon$
\begin{eqnarray}\label{1eq51}
(A^\varepsilon(t))^\dag && =
\varepsilon^2\int_0^{t/\varepsilon^2}a(\xi^{\varepsilon}(s),x^{\varepsilon}(s))d\nu(s)\leq
L\varepsilon^2\int_0^{t/\varepsilon^2}(1+(\xi^{\varepsilon}(s))^\dag)d\nu(s)\end{eqnarray}

Now, by Doob's inequality (see, e.g., \cite[Theorem 1.9.2]{lip1}),
$$\mathbf{E}((M_t^{\varepsilon})^\dag)^2\leq
4\abso{\mathbf{E}\langle M^{\varepsilon}\rangle_t},$$ (\ref{1qch})
and condition \textbf{C4} we obtain
\begin{eqnarray}\label{1eq6}
\abso{\langle
M^{\varepsilon}\rangle_t}=\abso{\varepsilon^2\int_0^{t/\varepsilon^2}\int_{\mathbb{R}^d\setminus
\{0\}}vv^*G(\xi^{\varepsilon}(s),dv;{x}_s^{\varepsilon})d\nu(s)}=\abso{\varepsilon^2\int_0^{t/\varepsilon^2}
c(\xi^{\varepsilon}(s);{x}^{\varepsilon}(s))d\nu(s)}\leq\nonumber\\
L\varepsilon^2\int_0^{t/\varepsilon^2}[1+((\xi^{\varepsilon}(s))^\dag)^2]d\nu(s).
\end{eqnarray}

Inequalities (\ref{1eq4})-(\ref{1eq6}), condition (\ref{1con1}) and
Cauchy-Bunyakovsky-Schwarz inequality,
([$\int_0^t\varphi(s)ds]^2\leq t\int_0^t\varphi^2(s)ds$), imply
$$\mathbf{E}((\xi^{\varepsilon}(t))^\dag)^2\leq
k_1+k_2\varepsilon^2\int_0^{t/\varepsilon^2}\mathbf{E}[((\xi^{\varepsilon}(s))^\dag)^2d\nu(s)]\leq
k_1+k_2l_1\varepsilon^2\int_0^{t/\varepsilon^2}\mathbf{E}((\xi^{\varepsilon}(s))^\dag)^2ds=$$
$$ k_1+k_2l_1\int_0^{t}\mathbf{E}((\xi^{\varepsilon}(s))^\dag)^2ds,$$
where $k_1, k_2$ and $l_1$ are positive constants independent of
$\varepsilon$.

By Gronwall inequality (see, e.g., \cite[p. 498]{ethier}), we obtain
$$\mathbf{E}((\xi^{\varepsilon}(t))^\dag)^2\leq k_1\exp(k_2l_1 t).$$

Thus, both parts of $\xi^{\varepsilon}(t)$ are compact and bounded,
so $$\mathbf{E}\sup\limits_{t\leq T}|\xi^{\varepsilon}(t)|^2\leq
k_T.$$

Hence the lemma is proved.
\fin

\begin{lm} Under assumption $\mathbf{C4}$ there exists a
constant $k>0$, independent of $\varepsilon$ such that
$$\mathbf{E}|\xi^{\varepsilon}(t)-\xi^{\varepsilon}(s)|^2\leq k |t-s|.$$
\end{lm}

\noindent{\it Proof}: In the same manner with (\ref{1eq4}), we may
write
$$|\xi^{\varepsilon}(t)-\xi^{\varepsilon}(s)|^2\leq 2|B_t^{\varepsilon}
-B_s^{\varepsilon}|^2+2|M_t^{\varepsilon}-M_s^{\varepsilon}|^2.$$ By
using Doob's inequality, we obtain
$$\mathbf{E}|\xi^{\varepsilon}(t)-\xi^{\varepsilon}(s)|^2\leq
2\mathbf{E}\{|B_t^{\varepsilon}-B_s^{\varepsilon}|^2+8\abso{\langle
M^{\varepsilon}\rangle_t-\langle M^{\varepsilon}\rangle_s}\}.$$

Now (\ref{1eq6}) and condition (\ref{1con1}) and assumption
$\mathbf{C4}$ imply
$$|B_t^{\varepsilon}-B_s^{\varepsilon}|^2+8\abso{\langle
M^{\varepsilon}\rangle_t-\langle M^{\varepsilon}\rangle_s}\leq
k_3[1+((\xi^{\varepsilon}(T))^\dag)^2]|t-s|,$$ where $k_3$ is a
positive constant independent of $\varepsilon$.

From the last inequality and Lemma 1 the desired conclusion is obtained.
\fin

The conditions proved in Corollary 2 and Lemma 2 are necessary and
sufficient for the compactness of the family of processes
$\xi^{\varepsilon}(t), t\geq 0, \varepsilon>0$.

\vsp \noindent{\sc Step 2}. At the next step of proof we apply the
problem of singular perturbation to the generator of the process
$\xi^{\varepsilon}(t).$ To do this, we mention the following
theorem.
$C^2_0(\mathbb{R}^d\times E)$ is the space of real-valued twice
continuously differentiable functions on the first argument, defined
on $\mathbb{R}^d\times E$ and vanishing at infinity, and
$C(\mathbb{R}^d\times E)$ is the space of real-valued continuous
bounded functions defined on $\mathbb{R}^d\times E$.

\begin{thm}(\cite[Theorem 6.3]{korlim}) Let the following conditions hold
for a family of Markov processes $\xi^{\varepsilon}(t), t\ge 0,
\varepsilon>0$:
\begin{description}
\item[CD1:] There exists a family of test functions
$\varphi^{\varepsilon}(u, x)$ in $C^2_0(\mathbb{R}^d\times E)$, such
that $$\lim\limits_{\varepsilon\to 0}\varphi^{\varepsilon}(u, x) =
\varphi(u),$$ uniformly on $u, x.$

\item[CD2:] The following convergence holds
$$\lim\limits_{\varepsilon\to
0}\mathbf{L}^{\varepsilon}\varphi^{\varepsilon}(u, x) =
\mathbf{L}\varphi(u),$$ uniformly on $u, x$. The family of functions
$\mathbf{L}^{\varepsilon}\varphi^{\varepsilon}, \varepsilon>0$ is
uniformly bounded, and $\mathbf{L}\varphi(u)$ and
$\mathbf{L}^{\varepsilon}\varphi^{\varepsilon}$ belong to
$C(\mathbb{R}^d\times E)$.

\item[CD3:] The quadratic characteristics of the
martingales that characterize a coupled Markov process
$\xi^{\varepsilon}(t), x^{\varepsilon}(t), t\geq0, \varepsilon>0$
have the representation $\left\langle \mu^{\varepsilon}\right\rangle_t = \int^t_0
\zeta^{\varepsilon}(s)ds,$ where the random functions
$\zeta^{\varepsilon}, \varepsilon> 0,$ satisfy the condition
$$\sup\limits_{0\leq s \leq T} \mathbf{E}|\zeta^{\varepsilon}(s)|\leq
c < +\infty.$$

\item[CD4:] The convergence of the initial values holds
and
$$\sup\limits_{\varepsilon>0}\mathbf{E}|\zeta^{\varepsilon}(0)|\leq C
< +\infty.$$
\end{description}

Then the weak convergence
$$\xi^{\varepsilon}(t)\Rightarrow \xi(t),\quad \varepsilon\to 0,$$
takes place.
\end{thm}

We consider the the extended Markov renewal process
\begin{eqnarray}\label{ren}\xi_n^{\varepsilon},{x}_n^{\varepsilon}, \tau_n^{\varepsilon}, n\ge
0,\end{eqnarray} where $x_n^{\varepsilon} =
x^{\varepsilon}(\tau_n^{\varepsilon}), x^{\varepsilon}(t) :=
x(t/\varepsilon^2),
\xi_n^{\varepsilon}=\xi^{\varepsilon}(\tau_n^{\varepsilon})$ and
$\tau_{n+1}^{\varepsilon} = \tau_{n}^{\varepsilon} +
\varepsilon^2\theta_n^{\varepsilon} , n \geq 0,$ and
$$P(\theta_{n+1}^{\varepsilon} \leq t | x_n^{\varepsilon} = x) =
F_x(t) = P(\theta_x \leq t).$$

\begin{dfn} \cite{svirid1} The \textit{compensating operator} $\mathbf{L}^{\varepsilon}$ of the Markov
renewal process (\ref{ren}) is defined by the following relation
$$\mathbf{L}^{\varepsilon}\varphi(\xi^{\varepsilon}_0,x_0,\tau_0) = q(x_0)\mathbf{E}[
\varphi(\xi^{\varepsilon}_1,x_1,\tau_1)
-\varphi(\xi^{\varepsilon}_0,x_0,\tau_0) | \mathcal{F}_0], $$ where
$$ \mathcal{F}_t :=
\sigma(\xi^{\varepsilon}(s),x^{\varepsilon}(s),\tau^{\varepsilon}(s);
0 \leq s \leq t).$$
\end{dfn}

Using Lemma 9.1 from \cite{korlim} we obtain that the compensating
operator of the extended Markov renewal process from Definition 1
can be defined by the relation (see also Section 2.8 in
\cite{korlim})
\begin{eqnarray}\label{1eq7}
\mathbf{L}^{\varepsilon}\varphi(u,v;x)=\varepsilon^{-2}q(x)\left[
\int_E P(x, dy) \int_{\mathbb{R}^d}
G_{u,x}^{\varepsilon}(dz)\varphi(u + z,v;y)-\right.\\\nonumber
\left.\varphi(u,v;x)\right].
\end{eqnarray}

By analogy with \cite[Lemma 9.2]{korlim} we may prove the following
result:

\begin{lm} The main part in the asymptotic representation of the compensating
operator (\ref{1eq7}) is as follows
$$\mathbf{L}^{\varepsilon}\varphi(u,v,x) = \varepsilon^{-2}\mathbf{Q}\varphi(\cdot,\cdot,x) +
\varepsilon^{-1}a_1(u;x)\mathbf{Q}_0\varphi'_u(u,\cdot,\cdot) +
[a(u;x) - a_0(u;x)]\mathbf{Q}_0\varphi'_u(u,\cdot,\cdot) +$$ $$
\frac{1}{2} [c(u;x)-c_0(u;x)]\mathbf{Q}_0\varphi_{uu}''(u,
\cdot,\cdot)+\mathbf{G}_{u,x}\mathbf{Q}_0\varphi(u,\cdot,\cdot)
$$
where:
$$\mathbf{Q}_0\varphi(x) := q(x) \int_E P(x, dy)\varphi(y),
\mathbf{G}_{u,x}\varphi(u) := \int_{\mathbb{R}^d} [\varphi(u + z)
-\varphi(u)]G_{u,x}(dz),$$ $$a_0(u;x)=\int_EvG_{u,x}(dv),
c_0(u;x)=\int_Evv^*G_{u,x}(dv).$$
\end{lm}

\noindent{\it Proof} of this Lemma is analogical to the proof of
\cite[Lemma 9.2]{korlim}.

The solution of the singular perturbation problem at the test
functions
$\varphi^{\varepsilon}(u,x)=\varphi(u)+\varepsilon\varphi_1(u,x)+\varepsilon^2\varphi_2(u,x)$
in the form
\begin{eqnarray}\label{spp}\mathbf{L}^{\varepsilon}\varphi^{\varepsilon}
={\mathbf{L}}\varphi+\theta^{\varepsilon}\varphi\end{eqnarray} can
be found in the same manner with Lemma 9.3 in \cite{korlim}.

To simplify the formula, we refer to the embedded Markov chain.
Corresponding generator $\widetilde{\mathbf{Q}}:=P-I,$ and the
potential operator satisfies the correlation
$\widetilde{R}_0(P-I)=\widetilde{\Pi}-I.$

From (\ref{spp}) we obtain
$$\widetilde{\mathbf{Q}}\varphi=0,$$
$$\widetilde{\mathbf{Q}}\varphi_1+\mathbf{A}_1(x)P\varphi=0,$$
$$\widetilde{\mathbf{Q}}\varphi_2+\mathbf{A}_1(x)P\varphi_1+(\mathbf{A}(x)+\mathbf{C}(x)+\mathbf{G}_{u,x})P\varphi=m(x){\mathbf{L}}\varphi,$$
where $$\mathbf{A}(x)\varphi(u):=[a(u;x)-a_0(u;x)]\varphi'(u),
\mathbf{A}_1(x)\varphi(u):=a_1(u;x)\varphi'(u),$$
$$\mathbf{C}(x):=\frac{1}{2} [c(u;x)-c_0(u;x)]\varphi_{uu}''(u).$$

From the second equation we obtain
$\varphi_1=\widetilde{R}_0\mathbf{A}_1(x)\varphi,$ and substituting
it into the last equation we have:
$$\widetilde{\mathbf{Q}}\varphi_2+\mathbf{A}_1(x)P\widetilde{R}_0\mathbf{A}_1(x)\varphi+(\mathbf{A}(x)+\mathbf{C}(x)+\mathbf{G}_{u,x})\varphi=m(x){\mathbf{L}}\varphi.$$

As soon as $P\widetilde{R}_0=\widetilde{R}_0+\widetilde{\Pi}-I$ we
finally obtain
\begin{eqnarray}\label{1eq8}
q^{-1}{\mathbf{L}}=\widetilde{\Pi}[(\mathbf{A}(x)+\mathbf{C}(x)+\mathbf{G}_{u,x})+\mathbf{A}_1(x)\widetilde{R}_0\mathbf{A}_1(x)-\mathbf{A}^2_1(x)]\widetilde{\Pi}.
\end{eqnarray}

Simple calculations give us (\ref{1limgen}) from (\ref{1eq8}).

Now Theorem 2 can be applied.

We see from (\ref{1eq7}) and (\ref{1eq8}) that the solution of
singular perturbation problem for
$\mathbf{L}^{\varepsilon}\varphi^{\varepsilon}(u,v;x)$ satisfies the
conditions \textbf{CD1, CD2}. Condition \textbf{CD3} of this theorem
implies that the quadratic characteristics of the martingale,
corresponding to a coupled Markov process, is relatively compact.
The same result follows from the CCC (see Corollary 2 and Lemma 2)
by \cite{jacod1}. Thus, the condition \textbf{CD3} follows from the
Corollary 2 and Lemma 2. Due to \textbf{L1} the condition
\textbf{CD4} is also satisfied. Thus, all the conditions of above
Theorem 2 are satisfied, so the weak convergence
$\xi^{\varepsilon}(t)\Rightarrow \xi^0(t)$ takes place.

Theorem 1 is proved.\fin

{\it Acknowledgements.} The authors thank University of Bielefeld
for hospitality and financial support by DFG project 436 UKR
113/94/07-09.\par

\vsp
\hrule

\end{document}